\documentclass[10pt,upref]{amsart}
\usepackage{latexsym,amsmath,amstext,mathrsfs,amssymb,amsthm,amscd,amsopn,verbatim,graphics}
\usepackage[all,2cell,dvips]{xy}
\UseAllTwocells
\SilentMatrices

\usepackage{times}
\numberwithin{equation}{section}

\newcommand{\norm}[1]{\lvert\!\lvert#1\rvert\!\rvert}

\DeclareMathOperator{\id}{id}

\DeclareMathOperator{\Gr}{Gr}

\newcommand{\dd}{\mathrm{d}}

\newcommand{\unint}{\mathrm{I}}

\newtheorem{Thm}{Theorem}[section]

\theoremstyle{remark}
\newtheorem{Rem}[Thm]{Remark}
\newtheorem*{Ack}{Acknowledgment}

\theoremstyle{definition}

\newtheorem{Def}[Thm]{Definition}
\newtheorem{Exa}[Thm]{Example}

\newtheorem*{Reg*}{Regularization procedure}

\newcommand{\braket}[2]{\left\langle{\,{#1}\,,\,{#2}\,}\right\rangle}

\newcommand{\bbR}{{\mathbb{R}}}

\newcommand{\calI}{\mathcal{I}}

\newcommand{\sfO}{{\mathsf{O}}}

\newcommand{\sfM}{{\mathsf{M}}}

\newcommand{\Inject}[2]{\calI\!\!\!\ \left(\bbR^{#1},\! \bbR^{#2}\right)}

\newcommand\qq{\rm}

\newcommand\jmp[1]{{\qq J.\ Math.\ Phys.\ \bf #1}}

\newcommand\np[1]{{\qq Nucl.\ Phys.\ \bf #1}}

\newcommand\jdg[1]{{\qq J.\ Diff.\ Geom.\ \bf #1}}

\begin{document} 

\title[On the space of injective\ldots]{On the space of injective linear maps from $\bbR^d$ into $\bbR^m$}

\author[C.~A.~Rossi]{Carlo~A.~Rossi}
\address{Dept.\ of Mathematics --- Technion ---  
32000 Haifa --- Israel}  
\email{crossi@techunix.technion.ac.il}
\thanks{C.~A.~Rossi acknowledges partial support from the Aly Kaufman Fellowship}

\begin{abstract}
In this short note, we investigate some features of the space $\Inject{d}{m}$ of linear injective maps from $\bbR^d$ into $\bbR^m$; in particular, we discuss in detail its relationship with the Stiefel manifold $V_{m,d}$, viewed, in this context, as the set of orthonormal systems of $d$ vectors in $\bbR^m$.
Finally, we show that the Stiefel manifold $V_{m,d}$ is a deformation retract of $\Inject{d}{m}$.
One possible application of this remarkable fact lies in the study of perturbative invariants of higher-dimensional (long) knots in $\bbR^m$: in fact, the existence of the aforementioned deformation retraction is the key tool for showing a vanishing lemma for configuration space integrals {\`a} la Bott--Taubes (see~\cite{BT} for the $3$-dimensional results and~\cite{CR1},~\cite{C} for a first glimpse into higher-dimensional knot invariants).
\end{abstract}

\maketitle


\section{Introduction}\label{sec-introduction}
Configuration space integrals {\`a} la Bott--Taubes are one of the main tools in the study of the cohomology of the space of knots in $\bbR^m$ (in particular, in the study of so-called {\em perturbative knot-invariants} for $m=3$)~\cite{BT} and~\cite{CCL}, and of the space of higher-dimensional knots~\cite{B}, i.e.\ imbeddings of $S^{m-2}$ into $\bbR^m$, for $m$ odd. 
See also~\cite{CR1} for a first glance towards a theory of perturbative knot invariants for higher-dimensional knots, generalizing the approach by Bott~\cite{B}.

The main feature of the compactifications of (relative) configuration spaces \`a la Fulton--MacPherson, appearing in the aforementioned configuration space integrals, is that the resulting compact spaces are manifolds with corners, hence admitting a natural stratification, corresponding intuitively to the ``speed of convergence'' of points in the interior to points on the boundary.
Since the proof that configuration space integrals produce genuine knot invariants or cohomology classes on the space of knots makes use of Stokes' Theorem, one is mainly interested in the codimension-$1$ boundary faces. 
Any such face corresponds to the collapse of $r\geq 2$ points together; boundary faces with only $2$ points collapsing together are ``nice'' and lead to a reformulation of the theory of knot invariants or cohomology classes of the space of knots in terms of a cohomology of (coloured) diagrams (see e.g.~\cite{BN} and~\cite{K}), whereas the remaining faces are ``bad''.
To deal with ``bad'' faces one is lead to formulate various {\em Vanishing Lemmata}: some of them are very easy to prove, and are essentially based on dimensional arguments, whereas others are more involved.
Typically, the latter vanishing lemmata are based, more or less directly, to the so-called {\em Kontsevich involution}~\cite{K}: e.g.\ in \cite{BT} and~\cite{CCL}, Kontsevich's involution is used directly to produce vanishing lemmata, while e.g.\ in~\cite{Th} and~\cite{CR2}, the authors use more involved involutions, steming from Kontsevich's one.
Still, as seen e.g.\ in~\cite{BT}, Kontsevich's involution is not always sufficient, since particularly bad faces can arise, e.g.\ when all points in the configuration space collapse together.
In such cases, one has to resort to other kinds of arguments in order to get rid of these bad faces: see e.g.\ again~\cite{BT} for more details about a possible way of circumvening such problems in the case of Chern--Simons invariants.

In~\cite{CR1}, a graph-cohomology for invariants of higher-dimensional knots was sketched thanks to new vanishing lemmata, based again on Kontsevich's involution; but we stumbled upon a particularly bad face, coming, as in $3$-dimensional Chern--Simons theory, from the collapse of all points in the compactified (relative) configuration spaces.
We mentioned the existence of a {\em deformations retraction} from the space of linear injective maps between Euclidean spaces to Stiefel manifolds as crucial ingredient in a possible ``cure'' for this boundary contribution, in the same spirit of~\cite{BT} for $3$-dimensional Chern--Simons theory.
The main purpose of this paper is to describe explicitly this deformation retraction.
The existence of such a deformation retract implies, among other things, that one can identify the corresponding de Rham cohomologies.
The main argument in the proof of the existence of the aforementioned deformation retraction is the {\em Gram--Schmidt orthonormalization procedure}: in fact, linear injective maps from $\bbR^d$ into $\bbR^m$ can be equivalently viewed as systems of $d$ {\em linearly independent} vectors in $\bbR^m$, hence such systems give rise, by the Gram--Schmidt procedure, to systems of $d$ orthonormal vectors in $\bbR^m$.
As a confirm that the Gram--Schmidt procedure gives rise to the correct deformation retraction from injective linear maps between Euclidean spaces and Stiefel manifolds, we recollect as a special the well-known fact that the $m$-dimensional sphere $S^m$ is a deformation retract of the Euclidean space $R^{m+1}$ without the origin.

\begin{Ack}
I cordially acknowledge Prof.\ A.\ S.\ Cattaneo for many different reasons: first of all, for having guided me attentively during my PhD years spent under his supervision, which brought me to the present result, and also for having read this manuscript carefully, pointing to me many small details, which I would have surely forgotten.
I also acknowledge the pleasant atmosphere of the Department of Mathematics of the Technion, Haifa, where this work was accomplished.
\end{Ack}

\section{The Gram--Schmidt map from $\Inject{d}m$ to the Stiefel manifold $V_{m,d}$}\label{sec-gramschmidt}
First of all, let us describe the main objects of this note, namely the Stiefel manifold $V_{m,d}$, for $d\leq m$, and the space of linear injective maps between Euclidean spaces $\bbR^d$ and $\bbR^m$, again for $d\leq m$.
\begin{Def}\label{def-stiefel}
The {\em Stiefel manifold $V_{m,d}$}, for any two positive integers $d$, $m$ satisfying $d\leq m$, is defined as the set of all {\em linear isometries of $\bbR^d$ into $\bbR^m$}.
\end{Def}
\begin{Rem}\label{rem-stiefelmetric}
One has to specify in advance Euclidean structures on both $\bbR^d$ and $\bbR^m$: we simply consider the standard Euclidean structures on both spaces.
\end{Rem}
Notice that, alternatively, $V_{m,d}$ can be viewed as the set of {\em orthonormal systems of $d$ vectors in $\bbR^m$}: for our purposes, we will mostly use this second characterization.
It can be proved, see e.g.~\cite{Hm}, that $V_{m,d}$ is indeed a smooth manifold.


\begin{Def}\label{def-inject}
The set of linear injective maps from $\bbR^d$ into $\bbR^m$, for any two integers $d$, $m$, obeying $d\leq m$, is denoted by $\Inject{d}{m}$.
\end{Def}
First of all, since $\bbR^d$, $\bbR^m$ are both endowed with their respective standard Euclidean structure, one can pick the standard orthonormal basis of $\bbR^d$, resp.\ $\bbR^m$, denoted by $\left\{e^d_1,\dots,e^d_d\right\}$, resp.\ $\left\{e^m_1,\dots,e^m_m\right\}$: hence, a general element $\alpha\in \Inject{d}{m}$ can be represented as a matrix with $d$ columns and $m$ rows as follows
\[
\alpha\equiv\left[\alpha\!(e_1^d),\dots,\alpha\!(e_d^d)\right],
\]
where $\alpha$ has been identified with the matrix whose column vectors $\alpha\!(e_i^d)$ with components
\[
\alpha\!(e_i^d)=\sum_{j=1}^d\alpha_{ij}e_j^m,\quad\forall 1\leq i\leq d.
\]
Thus, the column vectors of the matrix $\alpha$ are linearly independent vectors, represented as $m$-tuples with real coefficients w.r.t.\ the standard basis of $\bbR^m$. 
Hence, depending only upon the choice of two bases of $\bbR^d$ and $\bbR^m$ respectively, the space of linear injective maps $\Inject{d}{m}$ can be alternatively viewed also as the set of systems of $d$ linearly independent vectors in $\bbR^m$.
This implies, among othet things, that $\Inject{d}{m}$ is also a smooth manifold, being the complement of a closed subset of $\bbR^{d\times m}$, thus an open subset of $\bbR^{d\times m}$.

Since the Stiefel manifold $V_{m,d}$ can be alternatively seen as the set of orthonormal systems of $d$ vectors in $\bbR^m$, it can be thought of as a submanifold of $\Inject{d}{m}$ via the natural inclusion
\begin{equation}\label{eq-stiefelincl}
\iota_{m,d}\colon V_{m,d}\hookrightarrow \Inject{d}{m}.
\end{equation}

The purpose now is to construct explicitly a smooth map $\Lambda_{m,d}$ from $\Inject{d}{m}$ to $V_{m,d}$: such a map will be then proved to be a deformation retraction.
Thus, the Stiefel manifold $V_{m,d}$ is a deformation retract of $\Inject{d}m$.
The main ingredient to construct the map $\Lambda_{m,d}$ is the {\em Gram--Schmidt orthonormalization procedure}, which we will briefly recall prior to apply it to our situation.

The main ingredient in the Gram--Schmidt orthonormalization procedure (for some Euclidean vector space $V$) is a system of linearly independent vectors $\left\{v_1,\dots,v_n\right\}$ of $V$.
Notice that $V$ can be in principle a Hilbert space of infinite dimension.
The Gram-Schmidt orthogonalization procedure produces out of $\left\{v_1,\dots,v_n\right\}$, an orthonormal system $\left\{e_1^v,\dots,e_n^v\right\}$ via the inductive rule
\begin{equation}\label{eq-gramschmidt}
e_1^v\colon=\frac{v_1}{\norm{v_1}},\quad
\widetilde{e}_i^v\colon=v_i-\sum_{j=1}^{i-1}\braket{e_j^v}{v_i}e_j^v,\quad e_i^v\colon=\frac{\widetilde{e}_i^v}{\norm{\widetilde{e}_i^v}},
\end{equation} 
where by $\braket{\ }{\ }$, resp.\ $\norm{\ }$, we denoted the chosen Euclidean structure on $V$, resp.\ the induced norm.

We consider the particular case $V=\bbR^m$: $\braket{\ }{\ }$ is then the standard Euclidean scalar product on $\bbR^m$ and 
\[
\left\{v_1,\dots,v_d\right\}\colon=\left\{\alpha\!(e_1^d),\dots,\alpha\!(e_d^d)\right\},\quad \alpha\in \Inject{d}m.
\]
Accordingly, denote by
\begin{equation*}
\left\{e_1^\alpha,\dots,e_d^\alpha\right\}
\end{equation*}
the orthonormal set constructed via the Gram--Schmidt procedure (\ref{eq-gramschmidt}) out of the above system of linearly independent vectors. 
From the explicit formula (\ref{eq-gramschmidt}), one can see directly that the orthonormal system $\left\{e_1^{\alpha},\dots,e_d^\alpha\right\}$ depends {\em smoothly} on $\alpha\in\Inject{d}m$.
It makes thus sense to define the {\em Gram--Schmidt map} $\Lambda_{m,d}$ from $\Inject{d}m$ to $V_{m,d}$ as 
\begin{equation*}
\begin{aligned}
\Lambda_{m,d}\colon \Inject{d}{m}&\to V_{m,d}\\ 
\alpha\equiv\left\{\alpha\!(e_1^d),\dots,\alpha\!(e_d^d)\right\}&\mapsto \left\{e_1^\alpha,\dots,e_d^\alpha\right\}.
\end{aligned}
\end{equation*}
The following Theorem contains the main feature of the Gram--Schmidt map $\Lambda_{m,d}$. 
\begin{Thm}\label{thm-gramretract}
The Gram--Schmidt map $\Lambda_{m,d}$ is a deformation retraction from $V_{m,d}$ to $\Inject{d}m$: in other words, if we consider the natural inclusion (\ref{eq-stiefelincl}), there exists a smooth map $\widehat{\Lambda}_{m,d}$ from $\Inject{d}m\times\unint$ to $\Inject{d}m$, $\unint$ being from now on the unit interval $[0,1]$, satisfying 
\begin{equation*}
\begin{cases}
\widehat{\Lambda}_{m,d}(\alpha,0)&=\alpha,\\ 
\widehat{\Lambda}_{m,d}(\alpha,1)&=\left(\iota_{m,d}\circ\Lambda_{m,d}\right)\!(\alpha),\quad \forall \alpha\in \Inject{d}m.
\end{cases}
\end{equation*}
\end{Thm}
\begin{proof}
Our proof is subdivided in many steps.
\begin{itemize}
\item[{\bf(First step)}] We prove the following formula: for any chosen index $i\in\left\{1,\dots,d\right\}$, the unit vector $e_i^\alpha$ can be written as a linear combination of the vectors $\alpha\!(e_j^d)$ as follows:
\begin{equation}\label{eq-gramlincomb}
e_i^\alpha=\sum_{j=1}^i \lambda_{ij}\!(\alpha)\alpha\!(e_j^d),\quad\text{such that}\quad\lambda_{ii}\!(\alpha)>0.
\end{equation}
The coefficients in (\ref{eq-gramlincomb}) are labelled by $\alpha$, because they depend explicitly on it by (\ref{eq-gramschmidt}); their dependence on $\alpha$ is smooth, but not linear, as we will soon see.
Equation (\ref{eq-gramlincomb}) can be proved in an inductive way.
First of all, we consider $i=1$.
By (\ref{eq-gramschmidt}), one sees immediately that
\[
e_1^\alpha=\frac{\alpha\!(e_1^d)}{\norm{\alpha\!(e_1^d)}},\quad\text{thus it makes sense to set}\quad \lambda_{11}\!(\alpha)\colon=\frac{1}{\norm{\alpha\!(e_1^d)}}>0.
\] 
Clearly, the dependence on $\alpha$ is smooth.
Let us proceed with the inductive step: assume, for any $1\leq j\leq i$, (\ref{eq-gramlincomb}) holds, with the coefficients $\lambda_{jk}\!(\alpha)$ depending smoothly on $\alpha$.
Consider the explicit expression for $e_{i+1}^\alpha$:
\[
\widetilde{e}_{i+1}^\alpha\colon=\alpha\!(e_{i+1}^d)-\sum_{j=1}^i\braket{\alpha\!(e_{i+1}^d)}{e_j^\alpha}e_j^\alpha,\quad e_{i+1}^\alpha=\frac{\widetilde{e}_{i+1}^\alpha}{\norm{\widetilde{e}_{i+1}^\alpha}}.
\]
We need to perform some manipulations on the second term in the previous equation:
\begin{align*}
\widetilde{e}_{i+1}^\alpha&=\alpha\!(e_{i+1}^d)-\sum_{j=1}^i\braket{\alpha\!(e_{i+1}^d)}{e_j^\alpha}e_j^\alpha=\\
&\overset{\text{By induction}}=\alpha\!(e_{i+1}^d)-\sum_{j=1}^i\sum_{k=1}^j\lambda_{jk}\!(\alpha)\braket{\alpha\!(e_{i+1}^d)}{e_j^\alpha}\alpha\!(e_k^d)=\\
&=\alpha\!(e_{i+1}^d)-\sum_{j=1}^i\left[\sum_{k=j}^i\braket{\alpha\!(e_{i+1}^d)}{e_k^\alpha}\lambda_{kj}\!(\alpha)\right]\alpha\!(e_j^d)=\\
&=\sum_{j=1}^{i+1}\widetilde{\lambda}_{i+1,j}\!(\alpha)\alpha\!(e_j^d),\quad \widetilde{\lambda}_{i+1,i+1}\!(\alpha)=1.
\end{align*}
Dividing the explicit expression computed above inductively for $\widetilde{e}_{i+1}^\alpha$ by its norm, we get the desired result simply by putting
\[
\lambda_{i+1,j}\!(\alpha)\colon=\frac{\widetilde{\lambda}_{i+1,j}\!(\alpha)}{\norm{\widetilde{e}_{i+1}^\alpha}},\quad \lambda_{i+1,i+1}\!(\alpha)=\frac{1}{\norm{\widetilde{e}_{i+1}^\alpha}}>0.
\] 
By induction, it is also clear that the coefficients $\lambda_{i+1,j}\!(\alpha)$ depend smoothly on $\alpha$.
\item[{\bf(Second step)}]
We consider the matrix associated to $\alpha$, namely 
\[
\left[\alpha\!(e_1^d),\dots,\alpha\!(e_d^d)\right].
\]
Formula (\ref{eq-gramlincomb}) yields immediately
\begin{equation}\label{eq-grammatrix}
\begin{aligned}
\left[e_1^\alpha,\dots,e_d^\alpha\right]&=\left[\lambda_{11}\!(\alpha)\alpha\!(e_1^d),\dots,\sum_{i=1}^d\lambda_{d,i}\!(\alpha)\alpha\!(e_i^d)\right]=\\
&=\left[\alpha\!(e_1^d),\dots,\alpha\!(e_d^d)\right]\begin{bmatrix}
\lambda_{11}\!(\alpha)& \lambda_{21}\!(\alpha) & \cdots & \lambda_{d1}\!(\alpha)\\
0 & \lambda_{22}\!(\alpha) & \cdots & \lambda_{d2}\!(\alpha)\\
\vdots  & \cdots  & \ddots & \vdots \\
0 &\cdots &0& \lambda_{dd}\!(\alpha)
\end{bmatrix}.
\end{aligned}
\end{equation}
The matrix $\sfM(\alpha)$ on the right of $\left[\alpha\!(e_1^d),\dots,\alpha\!(e_d^d)\right]$ in the second equality in (\ref{eq-grammatrix}) is an upper triangular matrix, depending smoothly on $\alpha$ by our previous computations, whose eigenvalues are all strictly positive.
It is also clear that, if $\alpha$ is an isometry of $\bbR^d$ into $\bbR^m$ w.r.t.\ the standard Euclidean structures on $\bbR^d$ and $\bbR^m$ (hence, an element of $V_{m,d}$), then the Gram--Schmidt procedure leads to the matrix $\left[\alpha\!(e_1^d),\dots,\alpha\!(e_d^d)\right]$, whence it follows that in this case $\sfM(\alpha)=\id$.
Thus, the Gram--Schmidt procedure associates to a linear injective map $\alpha$ from $\bbR^d$ into $\bbR^m$ an upper triangular square $d\times d$-matrix $\sfM(\alpha)$ with strictly positive eigenvalues, which reduces to the identity when $\alpha$ is an isometry. 
\item[{\bf(Third step)}]
Using the map $\alpha\to\sfM(\alpha)$, we consider the family of upper triangular square $d\times d$-matrices associated to $\alpha$ in $\Inject{d}m$: 
\[
\sfM_t(\alpha)\colon=(1-t)\id+t\sfM(\alpha),\quad t\in\unint. 
\]  
Since $t\in\unint$, the matrix $\sfM_t(\alpha)$ has strictly positive eigenvalues, for any $t$ in $\unint$: in particular, it is a family of invertible $d\times d$-matrices, depending moreover smoothly on $t$.
We can finally consider the family of $m\times d$-matrices 
\begin{equation}\label{eq-gramfamily}
\widehat{\sfM}_t\!(\alpha)\colon=\left[\alpha\!(e_1^d),\dots,\alpha\!(e_d^d)\right]\left[(1-t)\id+t\sfM\!(\alpha)\right],\quad \alpha\in \Inject{d}m.
\end{equation}
For any $t\in\unint$, the column vectors are linearly independent (since the family $\sfM_t(\alpha)$ consists of invertible matrices): thus, for any $t\in \unint$, $\widehat{\sfM}_t(\alpha)$ can be viewed as an element of $\Inject{d}m$.
Moreover, the family given by (\ref{eq-gramfamily}) has the obvious properties:
\begin{align*}
\widehat{\sfM}_0\!(\alpha)&=\left[\alpha\!(e_1^d),\dots,\alpha\!(e_d^d)\right],\\ 
\widehat{\sfM}_1\!(\alpha)&=\left[\alpha\!(e_1^d),\dots,\alpha\!(e_d^d)\right]\sfM(\alpha)=\left[e_1^\alpha,\dots,e_d^\alpha\right].
\end{align*}
\item[{\bf(Fourth step)}] We finally define the map $\widehat{\Lambda}_{m,d}$ as follows:
\begin{equation}\label{eq-defretract}
\begin{aligned}
\widehat{\Lambda}_{m,d}:\Inject{d}m \times\unint&\to \Inject{d}m\\
(\alpha,t)\equiv(\left\{\alpha\!(e_1^d),\dots,\alpha\!(e_d^d)\right\},t)&\mapsto\left\{e_1^\alpha(t),\dots,e_d^\alpha(t)\right\},
\end{aligned}
\end{equation}
where $e_i^\alpha(t)$ denotes the $i$-th column vector of $\widehat{\sfM}_t(\alpha)$. 
By the previous computations, the map (\ref{eq-defretract}) satisfies
\begin{equation}\label{eq-propdefretr}
\widehat{\Lambda}_{m,d}\!(\alpha,0)=\left\{\alpha\!(e_1^d),\dots,\alpha\!(e_d^d)\right\}\equiv \alpha,\quad \widehat{\Lambda}_{m,d}\!(\alpha,1)=\left\{e_1^\alpha,\dots,e_d^\alpha\right\}=\left(\iota_{m,d}\circ\Lambda_{m,d}\right)\!(\alpha).
\end{equation}
Hence, the claim follows.
\end{itemize}
\end{proof}
Before ending this Section, let us illustrate Theorem~\ref{thm-gramretract} in two instructive special cases.
\begin{Exa}\label{exa-sphere}
Let us take a closer look at Theorem~\ref{thm-gramretract} in the particular case $d=1$ and $m\geq 1$.
First of all, we notice that the space $\Inject{1}{m+1}$ is obviously the Euclidean space $\bbR^{m+1}$ minus the origin, whereas the Stiefel manifold $V_{m+1,1}$ is the $m$-dimensional sphere $S^m$.
On the other hand, there is an obvious map from $\bbR^{m+1}\setminus\left\{0\right\}$ to $S^m$, the normalization map $v\overset{\rho}\mapsto v/\norm{v}$, for any nonzero $v$ in $\bbR^{m+1}$.
The normalization map corresponds to the easiest case of application of the Gram--Schmidt procedure.
Namely, the composition $\rho\circ\iota$ equals the identity map of $S^m$, while, on the other hand, the composite map $\iota\circ \rho$ takes the form
\[
\bbR^{m+1}\ni v\neq 0\mapsto\frac{v}{\norm{v}}\in S^m\subset\bbR^{m+1}\setminus\left\{0\right\}. 
\]
For any nonzero vector $v$ in $\bbR^{m+1}$, the family of $1\times1$-matrices (hence scalars) $\sfM_t$ takes the explicit form:
\[
\sfM_t(v)=(1-t)\id+t\sfM(v)=(1-t)+\frac{t}{\norm{v}}=t\left(\frac{1-\norm{v}}{\norm{v}}\right)+1.
\]
Therefore, the map $\widehat{\Lambda}_{m,1}$ takes the explicit form
\begin{equation*}
\widehat{\Lambda}_{m,1}(v)=\left[t\left(\frac{1-\norm{v}}{\norm{v}}\right)+1\right]v\in\Inject{1}{m}=\bbR^{m+1}\setminus\left\{0\right\},
\end{equation*} 
and this gives us the well-known fact that the $m$-dimensional sphere $S^m$ is a deformation retract of the Euclidean space $\bbR^{m+1}$ without the origin.

Hence, the Gram--Schmidt map $\Lambda_{m,d}$, for any $d\leq m$, from $\Inject{d}{m}$ to $V_{m,d}$ may be viewed as the natural generalization of the deformation retract from $\bbR^m\setminus\left\{0\right\}$ to the $m-1$-dimensional sphere $S^{m-1}$.
\end{Exa}

\begin{Exa}\label{exa-QRdecom}
Let us consider another special case, namely when $d=m$.
In this case, the space of injective linear maps $\Inject{m}m$ obviously equals the group $GL(m)$ of invertible real square $m\times m$-matrices; on the other hand, $V_{m,m}$ can be also obviously identified with the group $O(m)$ of orthogonal $m\times m$-matrices.
Given an element $\alpha$ of $GL(m)$, its column vectors form a basis of $\bbR^m$: the corresponding orthogonal matrix obtained by applying the Gram--Schmidt map $\Lambda_{m,m}$ has, as column vectors, the elements of the orthonormal basis obtained by the Gram--Schmidt procedure from the basis corresponding to the column vectors of $\alpha$; this orthogonal matrix we denote by $\Lambda_{m,m}(\alpha)$.
The key point in the proof of Theorem~\ref{thm-gramretract} is Formula (\ref{eq-grammatrix}): namely, in the case $d=m$, Formula (\ref{eq-grammatrix}) gives a decomposition of a general invertible $m\times m$-matrix $\alpha$ into an orthogonal $m\times m$-matrix $\Lambda_{m,m}(\alpha)$ and a regular upper triangular $m\times m$-matrix $\sfM(\alpha)$ with positive eigenvalues:
\[
\alpha=\Lambda_{m,m}(\alpha)\sfM(\alpha),\quad \alpha\in GL(m),\quad\Lambda_{m,m}(\alpha)\in O(m).
\]
But this is exactly the {\em QR decomposition} of an invertible matrix: thus, in the case $d=m$, the QR decomposition permits to construct explicitly a deformation retraction from $GL(m)$ to $O(m)$, as a corollary of Theorem~\ref{thm-gramretract}. 
\end{Exa}

\section{Equivariance properties of the Gram--Schmidt map}\label{sec-equivgramsch}
In this subsection, we briefly inspect equivariance properties of the map Gram--Schmidt map $\Lambda_{m,d}$ and of the corresponding homotopy $\widehat{\Lambda}_{m,d}$, introduced in Section~\ref{sec-gramschmidt}.

First of all, both manifolds $\Inject{d}{m}$ and $V_{m,d}$ carry, in a natural way, a left action of the group $SO(m)$ of orientation-preserving isometries of $\bbR^m$, given, respectively, by left multiplication of matrices, if we consider elements of $\Inject{d}{m}$ and $V_{m,d}$ as matrices, or by taking the diagonal action of the restriction of the standard action of $SO(m)$ on $\bbR^m$, if one looks at $\Inject{d}{m}$ and $V_{m,d}$ as $d$-tuples of $m$-dimensional vectors; it is obvious that both actions are compatible, hence one may switch from one action to the other without problems.
Obviously, the natural inclusion $\iota_{m,d}$ of the Stiefel manifold $V_{m,d}$ into $\Inject{d}m$ is $SO(m)$-equivariant.
The claim is that that both $\Lambda_{m,d}$ and the corresponding homotopy $\widehat{\Lambda}_{m,d}$ are also equivariant w.r.t.\ the aforementioned $SO(m)$-action, for any $d\leq m$.

Consider an object of $\Inject{d}m$ as a matrix $\alpha$ with corresponding system of $d$-linearly independent vectors in $\bbR^m$, $\left\{\alpha(e_1^d),\dots,\alpha(e_d^d)\right\}$.
Let $\sfO$ be a general element of $SO(m)$, whose action on $\Inject{d}m$ is given by left multiplication; the corresponding $d$-tuple $\left\{\alpha(e_1^d),\dots,\alpha(e_d^d)\right\}$ is therefore acted on by $\sfO$ as
\[
\left\{(\sfO\alpha)(e_1^d),\dots,(\sfO\alpha)(e_d^d)\right\}.
\]
We claim first that the following $SO(m)$-equivariance of the vectors $e_i^\alpha$ holds true:
\begin{equation}\label{eq-equivgramsch}
e_i^{\sfO\alpha}=\sfO(e_i^\alpha),\quad\forall \sfO\in SO(m),\quad \alpha\in \Inject{d}m,\quad \forall i=1,\dots,d.
\end{equation}
As in the first step of the proof of Theorem~\ref{thm-gramretract}, we use an inductive argument. 
The claim is obvious for $e_1^{\sfO\alpha}$, since $\sfO$ belongs to $SO(m)$, since we consider $\bbR^m$ endowed with the standard Euclidean metric.
Now, assume the $SO(m)$-equivariance $e_j^{\sfO\alpha}=\sfO(e_j^\alpha)$ holds true for all $j=1,\dots,i$, then one gets, by the properties of the Gram--Schmidt procedure,
\begin{align*}
\widetilde{e}_i^{\sfO\alpha}&=(\sfO\alpha)(e_i^d)-\sum_{j=1}^{i-1}\braket{(\sfO\alpha)(e_i^d)}{e_j^{\sfO\alpha}}e_j^{\sfO\alpha}=\\
&=\sfO(\alpha(e_i^d))-\sum_{j=1}^{i-1}\braket{\sfO(\alpha(e_i^d))}{\sfO(e_j^{\alpha})}\sfO(e_j^{\alpha})=\\
&=\sfO(\alpha(e_i^d))-\sum_{j=1}^{i-1}\braket{\alpha(e_i^d)}{e_j^{\alpha}}\sfO(e_j^{\alpha})=\\
&=\sfO\!\left(\widetilde{e}_i^\alpha\right),
\end{align*}
and since 
\[
e_i^\alpha=\frac{\widetilde{e}_i^\alpha}{\norm{\widetilde{e}_i^\alpha}},
\]
the claim follows immediately.
Consider now the claim proved in the first step of the proof of Theorem~\ref{thm-gramretract}, namely that, for any $\alpha$ in $\Inject{d}m$, there exist uniquely determined coefficients $\lambda_{ij}(\alpha)$, such that
\[
e_i^\alpha=\sum_{j=1}^i \lambda_{ij}(\alpha)\alpha(e_j^d),\quad i=1,\dots,m,\quad \lambda_{ii}(\alpha)>0.
\]
The coefficients $\lambda_{ij}(\alpha)$ are uniquely determined, since the system of vectors $\left\{e_1^\alpha,\dots,e_d^\alpha\right\}$ and $\left\{\alpha(e_1^d),\dots,\alpha(e_d^d)\right\}$ are both linearly independent.
Applying the previously showed $SO(m)$-invariance, we get
\[
e_i^{\sfO\alpha}=\sum_{j=1}^i \lambda_{ij}(\sfO\alpha)(\sfO\alpha)(e_j^d)\overset{!}=\sfO(e_i^\alpha)=\sum_{j=1}^i \lambda_{ij}(\alpha)(\sfO\alpha)(e_j^d),
\]
whence, since the coefficients $\lambda_{ij}$ are uniquely determined,
\[
\lambda_{ij}(\sfO\alpha)=\lambda_{ij}(\alpha),\quad \sfO\in SO(m),\quad \alpha\in\Inject{d}m.
\]
This implies immediately $SO(m)$-invariance of the matrix $\sfM(\alpha)$, introduced in the second step of the proof of Theorem~\ref{thm-gramretract}, whence it follows that the family $\sfM_t(\alpha)$ is also $SO(m)$-invariant.
Therefore, the consequence is that, by the $SO(m)$-invariance of $\sfM_t(\alpha)$, the map $\widehat{\Lambda}_{m,d}$, introduced in the fourth step of the proof of Theorem~\ref{thm-gramretract}, is $SO(m)$-equivariant.
All these computations can be summarized in the following
\begin{Thm}\label{thm-gramschequiv}
For any two integers $d\leq m$, the Gram--Schmidt map $\Lambda_{m,d}$ and the corresponding homotopy $\widehat{\Lambda}_{m,d}$ are both $SO(m)$-equivariant.
\end{Thm}

\section{Some consequences of the existence of the Gram--Schmidt map}\label{sec-cohoinject}
We come now to the main corollary of the existence of the Gram--Schmidt map $\Lambda_{m,d}$ and of the corresponding homotopy $\widehat{\Lambda}_{m,d}$, and the $SO(m)$-equivariance of both maps, namely that the ($SO(m)$-invariant) de Rham cohomologies of both spaces are isomorphic.
This is proved easily as follows: consider a differential form $\omega$ on $\Inject{d}m$ of degree $p$, then define 
\[
\widehat{\omega}_{m,d}\colon=\pi_{\unint*}\!\left(\widehat{\Lambda}_{m,d}^*(\omega)\right)\in \Omega^{p-1}\!\left(\Inject{d}m\right).
\]
In the previous formula, $\pi_{\unint*}$ denotes integration along the fiber of the (trivial) fibration $\Inject{D}m\times \unint\to \Inject{d}m$.
The (generalized) Stokes Theorem implies the following identity
\begin{equation}\label{eq-derhamgram}
\Lambda_{m,d}^*\!\left(\iota_{m,d}^*\omega\right)-\omega=\widehat{\dd \omega}_{m,d}+\dd \widehat{\omega}_{m,d},
\end{equation} 
whose immediate consequence is that the maps $\iota_{m,d}^*$ and $\Lambda_{d,m}$ are inverse one to another in de Rham cohomology.
Since one has~\ref{sec-equivgramsch} that both $\Inject{d}m$ and $V_{m,d}$ are left $SO(m)$-spaces, it makes sense to consider {\em left-invariant differential forms} on both spaces, i.e.\ differential forms $\omega$ on $\Inject{d}m$ or $V_{m,d}$ enjoying
\[
g^*\omega=\omega,\quad \forall g\in SO(m).
\]
Since the exterior derivative commutes with pull-backs, it makes sense to consider $SO(m)$-invariant de Rham cohomology of $\Inject{d}m$ and $V_{m,d}$: in this case, we consider de Rham cohomology classes of $SO(m)$-invariant forms.
(Caveat: notice that in $G$-invariant cohomology, a form is exact if and only if it is $G$-invariantly exact.)
Anyway, since $\iota_{m,d}$ and $\Lambda_{m,d}$ are $SO(m)$-equivariant, then they induce maps by pull-back on the invariant de Rham cohomology of the corresponding spaces.
Consider now a $SO(m)$-invariant form $\omega$ on $\Inject{d}m$ of degree $p$: then, the corresponding $p-1$-form $\widehat{\omega}_{m,d}$ is also $SO(m)$-invariant: this is a consequence of the fact that integration along the fiber of the trivial fibration $\Inject{d}m\times\unint\to \Inject{d}m$ commutes with the action of $SO(m)$, and commutes then with the homotopy $\widehat{\Lambda}_{m,d}$ by Theorem~\ref{thm-gramschequiv}.
This remarkable fact implies immediately that $\iota_{m,d}$ and $\Lambda_{m,d}$ induce (by pull-back) maps inverse one to another in $SO(m)$-invariant cohomology of the corresponding spaces.
Hence, the $SO(m)$-invariant de Rham cohomologies of $\Inject{d}m$ and $V_{m,d}$ are the same.
We are now going to discuss roughly, deserving to this topic a more detailed discussion somewhere else, the importance of this equality of (invariant) de Rham cohomologies.

Notice that the Stiefel manifold $V_{m,d}$ has a different characterization as a {\em homogeneous space}: in fact, it can be realized~\cite{Hm} as the quotient space $SO(m)/SO(m-d)$, $SO(m-d)$ imbedded as a subgroup of $SO(m)$ in the last $(m-d)\times (m-d)$ rows and columns.
(The characterization of $V_{m,d}$ as a homogeneous space gives a direct insight to the left action of $SO(m)$.)  
The $SO(m)$-invariant cohomology of $V_{m,d}$, see e.g.~\cite{GHV2}, can be computed easily using only algebraic tools, namely (relative) Lie algebra cohomology.
Thus, the a priori difficult to compute $SO(m)$-invariant cohomology of the space $\Inject{d}m$ can be computed easily: in particular, it is finite-dimensional.
In the following, we try and motivate the importance of spaces of injective linear maps and their invariant cohomologies.

Injective linear maps of $\bbR^d$ into $\bbR^m$ arise in the framework of (higher-dimensional) knot theory: namely, higher-dimensional knots can be viewed as imbeddings of a $d$-dimensional manifold $S$ (typically, $S$ can be the $d$-dimensional sphere $S^d$ or the Euclidean space $\bbR^d$) into a $m$-dimensional manifold $M$ (which is typically the Euclidean space $\bbR^m$), and, by their very definition, the tangent maps at any point of $S$ of such imbeddings are linear injective maps from the corresponding tangent spaces (which are non-canonically isomorphic to $\bbR^d$ and $\bbR^m$ respectively: they are in the case $S=\bbR^d$ and $M=\bbR^m$).
As a consequence, linear injective maps in $\Inject{d}m$ appear also in the compactification of relative configuration spaces $C_{p,q}\!\left(S,M\right)$, which are fibrations over the space of imbeddings of $S$ into $M$: the typical fibre over a given imbeddings consists of configuration of $p$ distinct points in $S$ and $q$ distinct points in $M$, such that the images of the points in $S$ w.r.t.\ the given imbedding are distinct from the points in $M$.
Such relative configuration spaces can be compactified \`a la Fulton--MacPherson: they were introduced and described in the framework of $3$-dimensional Chern--Simons theory in~\cite{BT} and also in~\cite{C} in the framework of higher-dimensional $BF$-theories.
The key point in the Bott--Fulton--MacPherson--Taubes compactification of configuration spaces is that the resulting compact manifolds are manifolds with corners: thus, they admit a natural stratification, corresponding, in rough terms, to ``rates of collapse'' of points.
Such rates of collapse are rigorously associated to limits; a detailed description of local coordinates of such compactifications may be found in~\cite{AS}, to which we refer for other details.

Now, the integrands involved in $3$-dimensional Chern--Simons theory and higher-dimensional $BF$-theories are products of so-called propagators, which are typically pull-backs of smooth forms on spheres (which are chosen to satisfy some invariance property, to which we will come later) w.r.t.\ maps of the form
\[
(x_i,x_j)\mapsto \frac{x_i-x_j}{\norm{x_i-x_j}},\quad x_i\neq x_j\in S\quad\text{or}\quad M.
\]
(Notice that one can take in principle also differences of a point in $S$ and of a point in $M$, or of two points in $S$ viewed as points in $M$, where any point in $S$ is mapped to a point in $M$ by some imbedding.)
This is the way differences of points on $S$ and on $M$ appear in configuration space integrals.
In order to show that the integrands are well-defined, i.e.\ they extend to smooth forms on compactified configuration spaces, one should look at the behavior of the integrands on the boundary, and this, in turn, is equivalent to show that the previous maps extend smoothly to the boundary.
To show this, we use the aforementioned limit procedures: in first approximation, differences between points, where an imbedding appears explicitly, are given by tangent maps at the point of collapse of imbeddings, thus to linear injective maps (roughly speaking, one approximates the imbedding by its Taylor expansion of order $1$).
 
Thus, applying the (generalized) Stokes Theorem to such integrals, and noting that the integrands are chosen to be closed forms, one is lead to the computation of integrals along boundary faces, and, by the rough arguments previously sketched, certain boundary faces can be viewed as pull-backs, w.r.t.\ ``tangent map'' of evaluation maps of imbeddings, of forms on the space of injective linear maps.
This ``tangent map'' obviously deals with tangent maps of imbeddings, i.e.\ with linear injective maps.

Usually, one has to consider explicitly such forms only in particular cases, since Vanishing Lemmata take care of the vanishing of most of such forms: however, as we sketched in some details in the final stages in~\cite{CR1} and in (still unpublished) details in~\cite{C}, there are certain forms on the space of linear injective maps which require a particular treatment.
This is the case, in Chern--Simons theory~\cite{BT}, of the most degenerate face, where {\em all points} collapse, and the same problem arises in higher-dimensional $BF$-theories~\cite{CR1} also when all points collapes together.
  
In~\cite{C} and~\cite{CR1}, we followed a strategy proposed by Bott and Taubes in Chern--Simons theory~\cite{BT}: namely, we showed first, by looking at boundary faces of codimension $2$, that such forms are closed.
Moreover, by their very construction, they are $SO(m)$-invariant; when restricted to Stiefel-manifolds, such forms are shown to be $SO(d)$-basic ($SO(d)$ acts on the right of $V_{m,d}=SO(m)/SO(m-d)$, and the quotient of this action can be identified with the Grassmann manifold $\Gr_{m,d}$), thus descending to $SO(m)$-invariant forms on Grassmann manifolds.
We mentioned, without explicitly giving the formula, the existence of the deformation retraction $\Lambda_{m,d}$: we use it to ``correct'' the aforementioned closed forms in the space of injective linear maps by the addition of exact forms, in order to get forms on Stiefel manifolds (this can be done using Formula (\ref{eq-derhamgram}).
As mentioned above, the ``corrected forms'' descend to $SO(m)$-invariant closed forms on Grassmann manifolds: such forms can be also explicitly computed by purely algebraic tools (Chevalley--Eilemberg Lie algebra cohomology).

We plan to write down all the details concerning the previous discussion elsewhere~\cite{CR2}, but the sketchy argument provided above should give the reader an idea of how important the $SO(m)$-invariant de Rham cohomology of spaces of linear injective maps is, and why the existence of a deformation retraction from such spaces to Stiefel manifolds plays a considerable r\^ole in the search for invariants of higher-dimensional knots and, more generally, cohomology classes on the spaces of knots.
On the other hand, we think that the result is interesting in itself, as it generalizes~\ref{exa-sphere} a well-known result for spheres, and produces an interesting application of the QR decomposition~\ref{exa-QRdecom} in topology. 

\end{document}